\documentclass[9pt]{article}
\usepackage{amsfonts,amssymb,amsmath,comment}
\usepackage{graphicx,graphics,color}

\usepackage{listings}

\usepackage{xcolor}
\definecolor{mygray}{rgb}{0.95,0.95,0.95}

 \makeatletter
 \@addtoreset{equation}{section}
 \makeatother

 \newcounter{enunciato}[section]

 \newtheorem{ittheorem}{Theorem}
 \newtheorem{itlemma}{Lemma}
 \newtheorem{itproposition}{Proposition}
 \newtheorem{itdefinition}{Definition}
 \newtheorem{itcorollary}{Corollary}

 \newtheorem{itconjecture}{Conjecture}

 \newenvironment{theorem}{\addtocounter{enunciato}{1}
 \begin{ittheorem}}{\end{ittheorem}}

\newcommand{\halmos}{\rule{1ex}{1.4ex}}

\def \ba {\begin{array}}
\def \ea {\end{array}}

\def \R {{\mathbb R}}

\def \N {{\mathbb N}}

\date{19 March 2019}
\begin{document}
\title{On a P\'olya functional for rhombi, isosceles triangles, and thinning convex sets.}

\author{M. van den Berg** , V. Ferone*, C. Nitsch*, C. Trombetti*
 \\
\\
**School of Mathematics, University of Bristol\\
University Walk, Bristol BS8 1TW, UK\\
\texttt{mamvdb@bristol.ac.uk}\\ \\
*Universit\`a degli Studi di Napoli Federico II\\
Via Cintia, Monte S. Angelo, I-80126 Napoli, Italy\\
\texttt{vincenzo.ferone@unina.it}\\
\texttt{c.nitsch@unina.it}\\
\texttt{cristina.trombetti@unina.it}}

\maketitle

\begin{abstract}
Let $\Omega$ be an open convex set in $\R^m$ with finite width, and with boundary $\partial \Omega$. Let $v_{\Omega}$ be the torsion function for $\Omega$, i.e. the solution of $-\Delta v=1, v|_{\partial\Omega}=0$. An upper bound is obtained for the product of $\Vert v_{\Omega}\Vert_{L^{\infty}(\Omega)}\lambda(\Omega)$, where $\lambda(\Omega)$ is the bottom of the spectrum of the Dirichlet Laplacian acting in $L^2(\Omega)$. The upper bound is sharp in the limit of a thinning sequence of convex sets.
For planar rhombi and isosceles triangles with area $1$, it is shown that $\Vert v_{\Omega}\Vert_{L^{1}(\Omega)}\lambda(\Omega)\ge \frac{\pi^2}{24}$, and that this bound is sharp.

\vskip 0.5truecm
\noindent
{\it AMS} 2000 {\it subject classifications.} 49J45, 49R05, 35P15, 47A75, 35J25.


\end{abstract}


\section{Introduction\label{sec1}}
\label{Introduction}
Let $\Omega$ be an open set in Euclidean space $\R^m$, and with boundary $\partial\Omega$. We denote the bottom of the spectrum of the Dirichlet Laplacian acting in $L^2(\Omega)$ by
\begin{equation*}
\lambda(\Omega)=\inf_{\varphi\in H_0^1(\Omega)\setminus\{0\}}\frac{\displaystyle\int_\Omega|D\varphi|^2}{\displaystyle\int_\Omega \varphi^2}.
\end{equation*}

It was shown in \cite{vdBC} and \cite{vdB} that if
\begin{equation*}
\lambda(\Omega)>0,
\end{equation*}
then the torsion function, i.e. the unique weak solution of
 $v_{\Omega}:\Omega\mapsto \R^+$,
\begin{equation*}
-\Delta v=1,\, v|_{\partial\Omega}=0,
\end{equation*}
satisfies
\begin{equation}\label{e4}
1\le \lambda(\Omega)M(\Omega)\le c_m,
\end{equation}
where
\begin{equation*}
M(\Omega)=\Vert v_{\Omega}\Vert_{L^{\infty}(\Omega)}.
\end{equation*}
In \cite{HV} it was shown that
\begin{equation*}
c_m\le \frac18(m+(5(4+\log 2))^{1/2}m^{1/2}+8)
\end{equation*}
The sharp constant in the right-hand side of \eqref{e4} is not known. However, for an open ball $B\subset\R^2$, an open square $S$, and an equilateral triangle $E$,
\begin{equation}\label{e7}
\lambda(B)M(B)<\lambda(S)M(S)<\lambda(E)M(E).
\end{equation}
The fact that $\lambda(B)M(B)<\lambda(E)M(E)$ was shown in \cite{HLP}. The full inequality \eqref{e7} follows from numerical evaluation of the series for the square, pp. 275-277 in \cite{TG}.

In \cite{MvdB} it was shown that the left-hand side of \eqref{e4} is sharp: for $\epsilon>0, m\ge 2,$ there exists an open bounded and connected set $\Omega_{\epsilon}\subset\R^m$ such that
\begin{equation*}
\lambda(\Omega_{\epsilon})M(\Omega_{\epsilon})<1+\epsilon.
\end{equation*}
For open, bounded convex sets in $\R^m$ it was shown in (3.12) of \cite{P} that
\begin{equation}\label{e9}
\lambda(\Omega)M(\Omega)\ge \frac{\pi^2}{8},
\end{equation}
with equality in the limit of an infinite slab (the open set with finite width bounded by two parallel $(m-1)$-dimensional planes).
The latter assertion has been made precise in \cite{HLP} where it was shown that if
\begin{equation}\label{e10}
S_n=(-n,n)^{m-1}\times(0,1),\, n\ge 1,
\end{equation}
then
\begin{equation}\label{e11}
\lambda(S_n)M(S_n)\le \frac{\pi^2}{8}+\frac{m-1}{8(n-\frac23)}.
\end{equation}
For bounded planar convex sets with  width $w(\Omega)$, and diameter $\textup{diam}(\Omega)$, it was shown in \cite{MvdB} that
\begin{equation}\label{e12}
\lambda(\Omega)M(\Omega)\le \frac{\pi^2}{8}\left(1+3^{2/3}7\left(\frac{w(\Omega)}{\textup{diam}(\Omega)}\right)^{2/3}\right).
\end{equation}

In Theorem \ref{the1} below we put \eqref{e10}-\eqref{e11}, and \eqref{e12} in a more general setting.
It was shown in Theorem 1.5 in \cite{vdBFNT} that for an open, bounded, convex set $\Omega$ with finite width $w(\Omega)$, and boundary $\partial\Omega$ there exist two points $z_0\in \partial \Omega$ and $z_1\in \partial\Omega$ such that $|z_0-z_1|=w(\Omega)$,
and such that the two hyper-planes tangent to $\partial\Omega$ through $z_0$ and $z_1$ are parallel. Denote these two hyper-planes by $H_0$ and $H_1$ respectively. Denote the inradius of $\Omega$, and the centre of an inball
 by $r(\Omega)$ and $c(\Omega)$, respectively:
\begin{equation*}
r(\Omega)=\sup_{x\in \Omega}\textup{dist}(x,\partial\Omega)=\textup{dist}(c(\Omega),\partial\Omega).
\end{equation*} We introduce Cartesian coordinates $(x',x_m) \in \R^m$ such that $x'(c(\Omega))=0, x_m(H_0)=0, x_m(H_1)=w(\Omega)$. Let  $H_{\mu}=\{x\in \R^m:x_m=\mu w(\Omega)\},\, 0\le\mu\le 1$, and let
\begin{equation*}
\Omega_{c}=\Omega\cap H_{x_m(c(\Omega))},
\end{equation*}
be the intersection of $\Omega$, and the hyper-plane through the centre of the inball and parallel to $H_0$.
We denote the inradius of this $(m-1)$-dimensional set by $\rho(\Omega)$. If the points $z_0,z_1,c(\Omega)$ are not unique then we choose them such that $\rho(\Omega)$ is maximal. The measure of $\Omega$ is denoted by $|\Omega|$.
\begin{theorem}\label{the1}
If $\Omega$ is an open, bounded, convex set in $\R^m, m\ge 2,$ then
\begin{equation}\label{e15}
\lambda(\Omega)M(\Omega)\le \frac{\pi^2}{8}\bigg(1+d_m\bigg(\frac{w(\Omega)}{\rho(\Omega)}\bigg)^{2/3}\bigg),
\end{equation}
where
\begin{equation}\label{e16}
d_m=  7(m+1)^{4/3}\pi^{-2}j_{(m-3)/2}^2,
\end{equation}
and where $j_{\nu}$ is the first positive zero of the Bessel function $J_{\nu}$.
\end{theorem}

The remaining results of this paper are for the P\'{o}lya functional for isosceles triangles and rhombi. Recall that the torsional rigidity (or torsion) $T(\Omega)$ of an open set $\Omega$ is defined by
\begin{equation}\label{e25}
T(\Omega)=\Vert v_{\Omega}\Vert_{L^1(\Omega)}=\int_{\Omega}v_{\Omega}.
\end{equation}
In P\'olya and Szeg\"o \cite{PSZ}, it was shown that for sets $\Omega$ with finite measure $|\Omega|$,
\begin{equation}\label{e26}
\frac{T(\Omega)\lambda(\Omega)}{|\Omega|}\le 1.
\end{equation}
The left-hand side of \eqref{e26} is the P\'olya functional for $\Omega$.
It was subsequently shown in \cite{vdBFNT} that the constant $1$ in the right-hand side above is sharp: for $\epsilon>0$, there exists an open, bounded, and connected set $\Omega_{\epsilon}\subset \R^m$ such that $T(\Omega_{\epsilon})\lambda(\Omega_{\epsilon})|\Omega_{\epsilon}|^{-1}\ge 1-\epsilon.$

The left-hand side of \eqref{e26} is invariant under the homothety transformation $t\mapsto t\Omega$. This implies for example that in Theorems \ref{the2}-\ref{the5} below we do not have to specify
the actual lengths of the edges of the rhombi and triangles. In the proofs of these theorems we fix the various lengths as a matter of convenience.

It was shown in Theorem 1.5 of \cite{vdBFNT} that for a thinning (collapsing) sequence $(\Omega_n)$ of bounded convex sets
\begin{equation}\label{e27}
\limsup_{n\rightarrow\infty}\frac{T(\Omega_n)\lambda(\Omega_n)}{|\Omega_n|}\le \frac{\pi^2}{12}.
\end{equation}
This supports the conjecture that for bounded, convex sets the sharp constant in the right-hand side of \eqref{e26} is $\pi^2/12$.

It was shown in Theorem 1.4 in \cite{vdBFNT} that for bounded convex sets in $\R^m, m\ge 3,$
\begin{equation}\label{e27a}
\frac{T(\Omega)\lambda(\Omega)}{|\Omega|}\ge \frac{\pi^2}{4m^{m+2}(m+2)},
\end{equation}
and that for planar, bounded, convex sets,
\begin{equation}\label{e27b}
\frac{T(\Omega)\lambda(\Omega)}{|\Omega|}\ge \frac{\pi^2}{48}.
\end{equation}

In Theorems \ref{the2}, \ref{the3}, \ref{the4}, and \ref{the5} we show that for isosceles triangles and rhombi the constant in the right-hand side of \eqref{e27b} can be improved to $\pi^2/24$, and that this constant is sharp.
\begin{theorem}\label{the2}
If $\triangle_{\beta}$ is an isosceles triangle with angles $\beta, \beta,\pi-2\beta$, and if $0<\beta\le \frac{\pi}{3}$ then
\begin{equation}\label{e28}
\frac{T(\triangle_{\beta})\lambda(\triangle_{\beta})}{|\triangle_{\beta}|}\le \frac{\pi^2}{24}\big(1+81 \big(\tan \beta\big)^{2/3}\big).
\end{equation}
\end{theorem}
\begin{theorem}\label{the3}
If $\lozenge_{\beta}$ is a rhombus with angles $\beta, \pi-\beta,\beta,\pi-\beta$, and if $\beta\le \frac{\pi}{3}$ then
\begin{equation}\label{e28a}
\frac{T(\lozenge_{\beta})\lambda(\lozenge_{\beta})}{|\lozenge_{\beta}|}\le\frac{\pi^2}{24}\big(1+15 \big(\tan\beta\big)^{2/3}\big).
\end{equation}
\end{theorem}
\begin{theorem}\label{the4}
If $\lozenge_{\beta}$ is as in Theorem \ref{the3}, then
\begin{equation}\label{e28b}
\frac{T(\lozenge_{\beta})\lambda(\lozenge_{\beta})}{|\lozenge_{\beta}|}\ge\frac{\pi^2}{24}.
\end{equation}
\end{theorem}
\begin{theorem}\label{the5}
If $\triangle_{\beta}$ is an isosceles triangle with angles $\beta,\beta, 2\pi-\beta$, then
\begin{equation}\label{e28i}
{ \frac{T(\triangle_{\beta})\lambda(\triangle_{\beta})}{|\triangle_{\beta}|}\ge \frac{\pi^2}{24}.}
\end{equation}
\end{theorem}

This paper is organised as follows. In Section \ref{sec2} we prove Theorem \ref{the1}. The proofs of Theorems \ref{the2} and \ref{the3} are deferred to Section \ref{sec3}. The proof of Theorem \ref{the4} is deferred to Section \ref{sec4}. { The proof of Theorem \ref{the5} consists of two parts. In Section \ref{sec5} part 1 we show that inequality \eqref{e28i} holds for all $\beta\in(0,\pi/3]\cup[\beta_0,\pi/2)$, where}
\begin{equation}\label{e28g}
\beta_0=\frac{\pi}{2}-\frac{33}{200}.
\end{equation}
In Section \ref{sec5} part 2 we use interval arithmetic to verify that \eqref{e28i} also holds for $\beta\in (\pi/3,\beta_0).$

\section{Proof of Theorem \ref{the1}\label{sec2}}
{\it Proof of Theorem} \ref{the1}.
We first observe, that by domain monotonicity of the torsion function, $v_{\Omega}$ is bounded by the torsion function for the (connected) set
bounded by $H_0$ and $H_1$. Hence
\begin{equation*}
v_{\Omega}(x)\le \frac12x_m(w(\Omega)-x_m)\le \frac{w(\Omega)^2}{8},\, (x',x_m)\in \Omega.
\end{equation*}
It suffices to obtain an upper bound for $\lambda(\Omega)$. By convexity we have that the convex hull of $z_0, z_1, \Omega_c$ is contained in $\Omega$. This convex hull in turn contains a cylinder with height $z\in [0,w(\Omega)]$, and base $\bigg(1-\frac{z}{w(\Omega)}\bigg)\Omega_c.$ Denote the first $(m-1)$-dimensional Dirichlet eigenvalue of $\Omega_c$ by $\lambda_c$. Then, by separation of variables, we have
\begin{equation}\label{e18}
\lambda(\Omega)\le \frac{\pi^2}{z^2}+\bigg(1-\frac{z}{w(\Omega)}\bigg)^{-2}\lambda_{c}.
\end{equation}
The right-hand side of \eqref{e18} is minimised for
\begin{equation*}
\frac{1}{z}=\frac{1}{w(\Omega)}+\bigg(\frac{\lambda_{c}}{\pi^2w(\Omega)}\bigg)^{1/3}.
\end{equation*}
This gives that
\begin{equation}\label{e20}
\lambda(\Omega)\le \frac{\pi^2}{w(\Omega)^2}\bigg(1+3\bigg(\frac{\lambda_{c}w(\Omega)^2}{\pi^2}\bigg)^{1/3}+
3\bigg(\frac{\lambda_{c}w(\Omega)^2}{\pi^2}\bigg)^{2/3}+\frac{\lambda_{c}w(\Omega)^2}{\pi^2}\bigg).
\end{equation}

The inball intersects $\Omega_c$ in a $(m-1)$-dimensional disc with radius $r(\Omega)$ which is, by a generalisation of Blaschke's theorem (see p.215 in \cite{Y}, and p.79 in \cite{E}), bounded from below by $w(\Omega)/(m+1)$.
Hence
\begin{equation}\label{e25}
\lambda_{c}\le (m+1)^2j_{(m-3)/2}^2w(\Omega)^{-2}.
\end{equation}
By \eqref{e20} and \eqref{e25} we obtain
\begin{align}\label{e26}
\lambda(\Omega)&\le \frac{\pi^2}{w(\Omega)^2}\bigg(1+\bigg(\frac{\lambda_{c}w(\Omega)^2}{\pi^2}\bigg)^{1/3}\bigg(3+3\bigg(\frac{(m+1)^2j_{(m-3)/2}^2}{\pi^2}\bigg)^{1/3}
+\bigg(\frac{(m+1)^2j_{(m-3)/2}^2}{\pi^2}\bigg)^{2/3}\bigg)\bigg)\nonumber \\ & \le \frac{\pi^2}{w(\Omega)^2}\bigg(1+7 \bigg(\frac{(m+1)^2j_{(m-3)/2}^2}{\pi^2}\bigg)^{2/3}\bigg(\frac{\lambda_{c}w(\Omega)^2}{\pi^2}\bigg)^{1/3}\bigg).
\end{align}
Since  the $(m-1)$- dimensional set $\Omega_c$ contains a disc of radius $\rho(\Omega)$ we have
\begin{equation}\label{e24}
\lambda_{c}\le j_{(m-3)/2}^2\rho(\Omega)^{-2},
\end{equation}
and \eqref{e15}, \eqref{e16} follows by \eqref{e26} and \eqref{e24}.
\hspace*{\fill }$\square $

\section{Proofs of Theorem \ref{the2} and Theorem \ref{the3} \label{sec3}}
{\it Proof of Theorem} \ref{the2}.
Let $\triangle_{\beta}$ be an isosceles triangle with a base of length $2$ and width (height) of length $d$, and angles $\beta,\beta$, and $\pi-2\beta$ respectively. By hypothesis, $\beta=\arctan d\le \pi/3$ so that
$d\le \sqrt3$. We denote the infinite sector with opening angle $\beta$ by
\begin{equation*}
\Omega_{\beta}=\{(r,\phi): r>0, -\beta/2<\phi<\beta/2\}.
\end{equation*}
It is straightforward to verify that the torsion function for $\Omega_{\beta}$ is given by,
\begin{equation*}
v_{\Omega_{\beta}}(r,\phi)=\frac{r^2}{4}\bigg(\frac{\cos(2\phi)}{\cos \beta}-1\bigg),\, r>0, -\beta/2<\phi<\beta/2.
\end{equation*}
Let
\begin{equation*}
R=\big(1+d^2\big)^{1/2}.
\end{equation*}
We can cover $\triangle_{\beta}$ with two sectors of opening angles $\beta$ and radii $R$ each. By monotonicity and positivity of the torsion function we have
\begin{align}\label{e31}
T(\triangle_{\beta})&=\int_{\triangle_{\beta}}v_{\triangle_{\beta}}\nonumber \\ &\le 2\int_0^R dr\, r \int_{-\beta/2}^{\beta/2} d\phi\, v_{\Omega_{\beta}}(r,\phi)\nonumber \\ &=\frac18\big(1+d^2\big)^2\big(\tan \beta-\beta\big)\nonumber \\ &=\frac18\big(1+d^2\big)^2(d-\arctan d)\nonumber \\ &\le
\frac{d^3}{24}\big(1+d^2\big)^2,
\end{align}
where we have used that $d-\arctan d\le d^3/3$.
By adapting formula (31) in the proof of Theorem 2 in \cite{MvdB} to the geometry of $\triangle_{\beta}$ we find that
\begin{equation}\label{e33}
\lambda(\triangle_{\beta})\le \frac{\pi^2}{d^2}\bigg(1+7\bigg(\frac{d}{2}\bigg)^{2/3}\bigg).
\end{equation}
By \eqref{e31}, \eqref{e33}, and $|\triangle_{\beta}|=d,$ we obtain
\begin{align*}
\frac{T(\triangle_{\beta})\lambda(\triangle_{\beta})}{|\triangle_{\beta}|}&\le\frac{\pi^2}{24}\big(1+d^2\big)^2\bigg(1+7\bigg(\frac{d}{2}\bigg)^{2/3}\bigg)\nonumber \\ &
\le \frac{\pi^2}{24}\big(1+81d^{2/3}\big)\nonumber \\ &
=\frac{\pi^2}{24}\big(1+81\big(\tan \beta)^{2/3}\big),\, 0<\beta\le \frac{\pi}{3}.
\end{align*}
\hspace*{\fill }$\square $

{\it Proof of Theorem} \ref{the3}. Let $\lozenge_{\beta}$ be a rhombus with angles $\beta,\pi-\beta,\beta,\pi-\beta$, and diagonals of length $2$ and $d$ respectively. By hypothesis we have that $\beta\le \pi/3$, and $d\le 2/\sqrt3.$ This rhombus is covered by two sectors of opening angle $\beta=2\arctan(d/2)$, and radius $R=(1+(d/2)^2)^{1/2}.$ By the calculations in the proof of Theorem \ref{the2} we find that
\begin{align*}
T(\lozenge_{\beta})&\le \frac18R^4 (\tan \beta-\beta)\nonumber \\ & =\frac18\bigg(1+\frac{d^2}{4}\bigg)^2\bigg(\frac{d}{1-\frac{d^2}{4}}-2\arctan \big(\frac{d}{2}\big)\bigg)\nonumber \\ &\le \frac18\bigg(1+\frac{d^2}{4}\bigg)^2\bigg(\frac{d}{1-\frac{d^2}{4}}-d+\frac{d^3}{12}\bigg)\nonumber \\ &\le
\frac{d^3}{24}\bigg(1+\frac{d^2}{4}\bigg)^2\bigg(1+\frac{9d^2}{32}\bigg),\, 0<d\le \frac23\sqrt3.
\end{align*}
By adapting formula (31) in the proof of Theorem 2 in \cite{MvdB} to the geometry of $\lozenge_{\beta}$ we find that
\begin{equation*}
\lambda(\lozenge_{\beta})\le \frac{\pi^2}{d^2}\bigg(1+7\bigg(\frac{d}{2}\bigg)^{2/3}\bigg).
\end{equation*}
This, together with $|\lozenge_{\beta}|=d$ gives that,
\begin{align*}
\frac{T(\lozenge_{\beta})\lambda(\lozenge_{\beta})}{|\lozenge_{\beta}|}&\le \frac{\pi^2}{24}\bigg(1+\frac{d^2}{4}\bigg)^2\bigg(1+\frac{9d^2}{32}\bigg)\bigg(1+7\bigg(\frac{d}{2}\bigg)^{2/3}\bigg)\nonumber \\ &
\le\frac{\pi^2}{24}\bigg(1+15\bigg(\frac{d}{2}\bigg)^{2/3}\bigg)\nonumber \\ &
=\frac{\pi^2}{24}\bigg(1+15\big(\tan\beta\big)^{2/3}\bigg),\, 0<\beta\le \frac{\pi}{3}.
\end{align*}
\hspace*{\fill }$\square $

\section{Proof of Theorem \ref{the4}\label{sec4}}

Let $\lozenge_{\beta}$ be a rhombus such that major and minor diagonals have lengths $2$ and $d$, respectively (see Figure \ref{rh1}). We want to estimate the torsion and to this aim we use a test function

\def\svgwidth{9. cm}
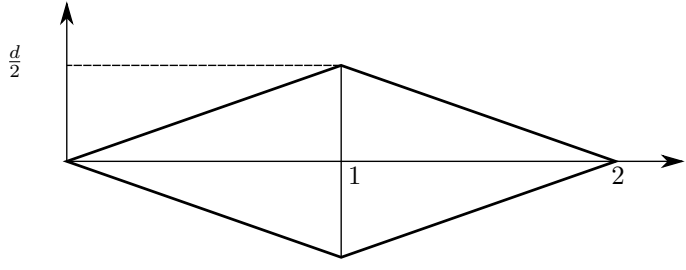
\begin{figure}
\centering
\begingroup%
  \makeatletter%
  \providecommand\color[2][]{%
    \errmessage{(Inkscape) Color is used for the text in Inkscape, but the package 'color.sty' is not loaded}%
    \renewcommand\color[2][]{}%
  }%
  \providecommand\transparent[1]{%
    \errmessage{(Inkscape) Transparency is used (non-zero) for the text in Inkscape, but the package 'transparent.sty' is not loaded}%
    \renewcommand\transparent[1]{}%
  }%
  \providecommand\rotatebox[2]{#2}%
  \ifx\svgwidth\undefined%
    \setlength{\unitlength}{394.74563478bp}%
    \ifx\svgscale\undefined%
      \relax%
    \else%
      \setlength{\unitlength}{\unitlength * \real{\svgscale}}%
    \fi%
  \else%
    \setlength{\unitlength}{\svgwidth}%
  \fi%
  \global\let\svgwidth\undefined%
  \global\let\svgscale\undefined%
  \makeatother%
  \begin{picture}(1,0.38076066)%
    \put(0,0){\includegraphics[width=\unitlength,page=1]{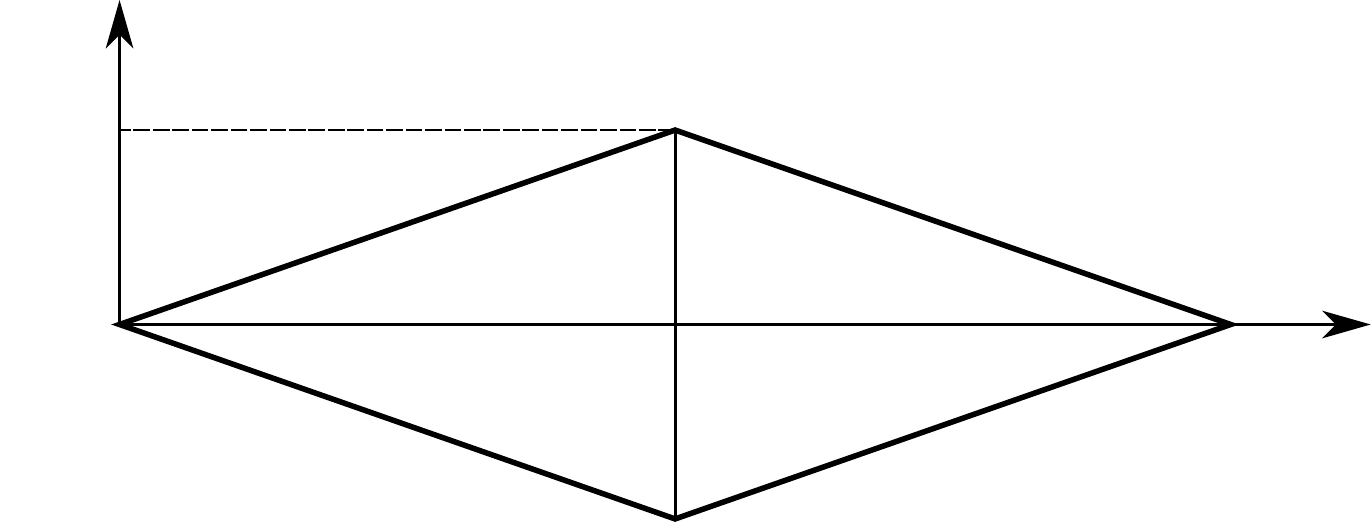}}%
    \put(-0.00252338,0.28153143){\color[rgb]{0,0,0}\makebox(0,0)[lb]{\smash{$\frac{d}{2}$}}}%
    \put(0.50268441,0.11071614){\color[rgb]{0,0,0}\makebox(0,0)[lb]{\smash{$1$}}}%
    \put(0.89208525,0.11071614){\color[rgb]{0,0,0}\makebox(0,0)[lb]{\smash{$2$}}}%
  \end{picture}%
\endgroup%

\caption{The rhombus of diagonals $2$ and $d$. To estimate the torsion we construct a test function $v$ which is symmetric with respect to the minor diagonal.}\label{rh1}
\end{figure}

\begin{equation*}
v(x,y) =
\left\{
\begin{array}{ll}
\displaystyle \frac{d^2 x^2}{4} - y^2,  & \mbox{ $0\le x\le 1$}, \\\\
\displaystyle \frac{d^2 (2-x)^2}{4} - y^2, & \mbox{ $1\le x\le 2$}.
\end{array}
\right.
\end{equation*}

In view of the variational definition of the torsion we have
\begin{equation*}
\frac{1}{T(\lozenge_{\beta})} \leq \frac{\displaystyle\int_{\lozenge_{\beta}}|Dv|^2}{\left(\displaystyle\int_{\lozenge_{\beta}}v\right)^2}=\frac{24+18d^2}{d^3}.
\end{equation*}

On the other hand we can estimate from below the first Dirichlet Laplacian eigenvalue of any rhombus by means of the Dirichlet Laplacian eigenvalue of a rectangle obtained by Steiner symmetrising the rhombus along a direction parallel to one of the sides (see Figure \ref{rh2}).
We denote by $b$ and $h$ the base and the height of the rectangle, respectively. Since the base $b$ coincides with the side of the rhombus, $b^2=1+{d^2\over 4}$, and $h=\frac{d}{\sqrt{1+{d^2\over 4}}}$.

We have,

\begin{equation*}
\lambda(\lozenge_{\beta}) \ge \pi^2 \left(\frac{1}{b^2}+\frac{1}{h^2}\right)=\pi^2 \frac{16+24d^2+d^4}{d^2(16+4d^2)}.
\end{equation*}
Observing that the area of the rhombus is equal to $d$, we have

\begin{equation}
\dfrac{\lambda(\lozenge_{\beta})T(\lozenge_{\beta})}{|\lozenge_{\beta}|} \ge \frac{\pi^2}{24}\frac{16+24d^2+d^4}{(1+\frac{3}{4}d^2)(16+4d^2)}\ge \frac{\pi^2}{24}, \qquad \mbox{ $0\le d\le 2$}.
\end{equation}
{\hspace*{\fill }$\square $}

\def\svgwidth{9. cm}
\begin{figure}
\centering

\begingroup%
  \makeatletter%
  \providecommand\color[2][]{%
    \errmessage{(Inkscape) Color is used for the text in Inkscape, but the package 'color.sty' is not loaded}%
    \renewcommand\color[2][]{}%
  }%
  \providecommand\transparent[1]{%
    \errmessage{(Inkscape) Transparency is used (non-zero) for the text in Inkscape, but the package 'transparent.sty' is not loaded}%
    \renewcommand\transparent[1]{}%
  }%
  \providecommand\rotatebox[2]{#2}%
  \ifx\svgwidth\undefined%
    \setlength{\unitlength}{338.53585485bp}%
    \ifx\svgscale\undefined%
      \relax%
    \else%
      \setlength{\unitlength}{\unitlength * \real{\svgscale}}%
    \fi%
  \else%
    \setlength{\unitlength}{\svgwidth}%
  \fi%
  \global\let\svgwidth\undefined%
  \global\let\svgscale\undefined%
  \makeatother%
  \begin{picture}(1,0.48958178)%
    \put(0,0){\includegraphics[width=\unitlength,page=1]{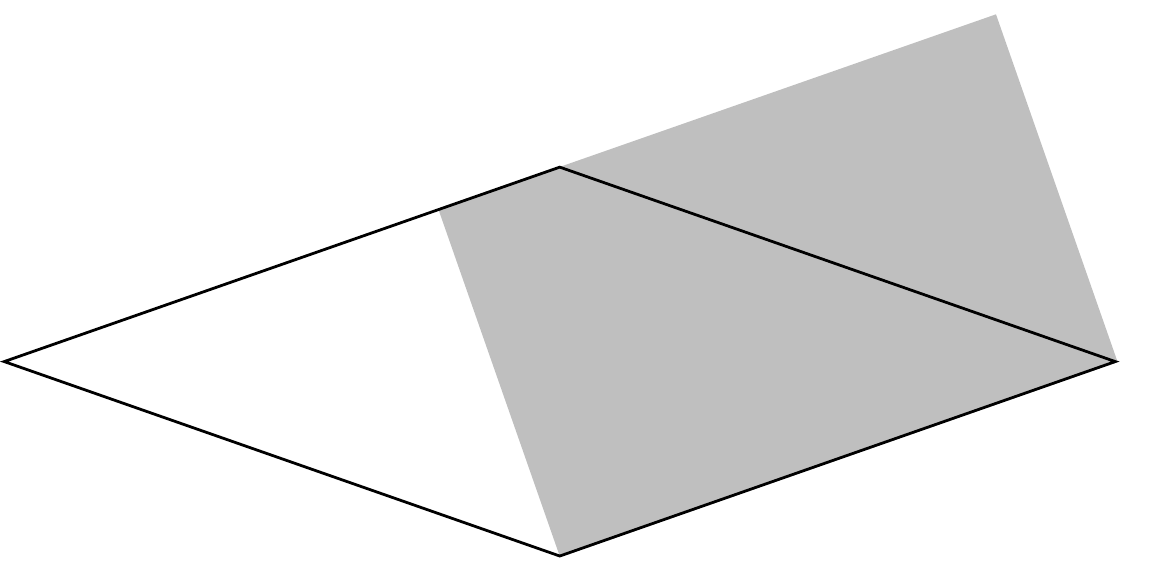}}%
    \put(0.73614309,0.05706348){\color[rgb]{0,0,0}\makebox(0,0)[lb]{\smash{$b$}}}%
    \put(0.93532004,0.34401357){\color[rgb]{0,0,0}\makebox(0,0)[lb]{\smash{$h$}}}%
    \put(0,0){\includegraphics[width=\unitlength,page=2]{rhombi2.pdf}}%
  \end{picture}%
\endgroup%

\caption{The rectangle shaded in grey is obtained by Steiner symmetrization. The Dirichlet Laplacian eigenvalue of the rectangle provides an estimate from below for the one on the rhombus.}\label{rh2}
\end{figure}

\section{Proof of Theorem \ref{the5}\label{sec5}}
\subsection{Proof for the case $\beta\in(0,\pi/3]\cup[\beta_0,\pi/2)$}
Let  $\triangle_{\beta}$ be an isosceles triangle with angles $\beta,\beta,\alpha=\pi-2\beta.$
We first consider the case $\dfrac{\pi}{3} \le \alpha< \pi$. We denote the height by $d$, and we fix the length of the basis equal to $2$. See Figure \ref{fig1}.

\def\svgwidth{9. cm}
\begin{figure}
\centering

\begingroup%
  \makeatletter%
  \providecommand\color[2][]{%
    \errmessage{(Inkscape) Color is used for the text in Inkscape, but the package 'color.sty' is not loaded}%
    \renewcommand\color[2][]{}%
  }%
  \providecommand\transparent[1]{%
    \errmessage{(Inkscape) Transparency is used (non-zero) for the text in Inkscape, but the package 'transparent.sty' is not loaded}%
    \renewcommand\transparent[1]{}%
  }%
  \providecommand\rotatebox[2]{#2}%
  \ifx\svgwidth\undefined%
    \setlength{\unitlength}{568.26546939bp}%
    \ifx\svgscale\undefined%
      \relax%
    \else%
      \setlength{\unitlength}{\unitlength * \real{\svgscale}}%
    \fi%
  \else%
    \setlength{\unitlength}{\svgwidth}%
  \fi%
  \global\let\svgwidth\undefined%
  \global\let\svgscale\undefined%
  \makeatother%
  \begin{picture}(1,0.39136836)%
    \put(0,0){\includegraphics[width=\unitlength,page=1]{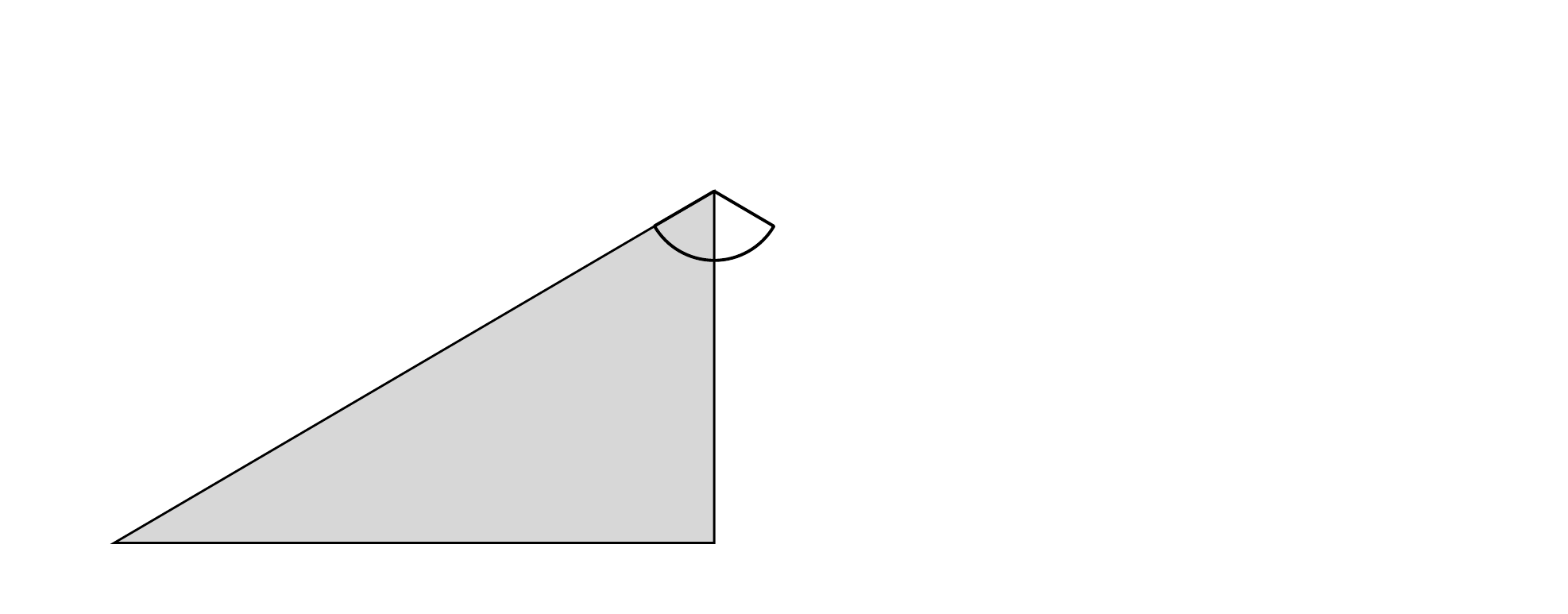}}%
    \put(-0.00321359,0.25859017){\color[rgb]{0,0,0}\makebox(0,0)[lb]{\smash{$h$}}}%
    \put(0.82524267,0.00702317){\color[rgb]{0,0,0}\makebox(0,0)[lb]{\smash{$2$}}}%
    \put(0,0){\includegraphics[width=\unitlength,page=2]{disegno1.pdf}}%
    \put(0.44394609,0.00568994){\color[rgb]{0,0,0}\makebox(0,0)[lb]{\smash{$1$}}}%
    \put(0.47007419,0.19539356){\color[rgb]{0,0,0}\makebox(0,0)[lb]{\smash{$\alpha$}}}%
    \put(0.17851772,-0.13975667){\color[rgb]{0,0,0}\makebox(0,0)[lb]{\smash{}}}%
    \put(0,0){\includegraphics[width=\unitlength,page=3]{disegno1.pdf}}%
  \end{picture}%
\endgroup%

\caption{Isosceles triangle of basis $2$, vertex angle $\alpha$ and height $h$. A test function to estimate the torsion is constructed on the shaded part and symmetrically reflected along the height.}\label{fig1}
\end{figure}
We use the function
\begin{equation*}
u(x,y) =
\left\{
\begin{array}{ll}
\displaystyle \frac{d^2 x^2}{4} - \left(y - \frac{dx}{2}\right)^2, & \mbox{ $0\le x\le 1$}, \\\\
\displaystyle \frac{d^2 (2-x)^2}{4} - \left(y - \frac{d(2-x)}{2}\right)^2, & \mbox{ $1\le x\le 2$},
\end{array}
\right.
\end{equation*}
as a test function for the torsion of $\triangle_{\beta}$. We find that
\begin{equation}
\label{T1}
\frac{2}{T(\triangle_{\beta})} \leq \frac{48(1+d^2)}{d^3}.
\end{equation}
Hence
\begin{equation}
\label{T2}
\frac{T(\triangle_{\beta})}{|\triangle_{\beta}|}\ge \frac{1}{24}\left(1 + \frac{1}{d^2}\right)^{-1} .
\end{equation}

We wish to estimate $\lambda(\triangle_{\beta})$ from below. To this aim we consider the first Dirichlet eigenfunction of $\triangle_{\beta}$ restricted to $x\in[0,1]$ and we reflect it, anti-symmetrically, with respect to the line $y=d x$ (see Figure \ref{fig2}).
This new function is a test function defined on the rectangle of sides $1,d$ (shaded in grey in Figure \ref{fig2}) orthogonal to the first eigenfunction of the Laplacian with the mixed boundary conditions described in Figure \ref{fig2}.

\def\svgwidth{9. cm}
\begin{figure}
\centering

\begingroup%
  \makeatletter%
  \providecommand\color[2][]{%
    \errmessage{(Inkscape) Color is used for the text in Inkscape, but the package 'color.sty' is not loaded}%
    \renewcommand\color[2][]{}%
  }%
  \providecommand\transparent[1]{%
    \errmessage{(Inkscape) Transparency is used (non-zero) for the text in Inkscape, but the package 'transparent.sty' is not loaded}%
    \renewcommand\transparent[1]{}%
  }%
  \providecommand\rotatebox[2]{#2}%
  \ifx\svgwidth\undefined%
    \setlength{\unitlength}{480.1634233bp}%
    \ifx\svgscale\undefined%
      \relax%
    \else%
      \setlength{\unitlength}{\unitlength * \real{\svgscale}}%
    \fi%
  \else%
    \setlength{\unitlength}{\svgwidth}%
  \fi%
  \global\let\svgwidth\undefined%
  \global\let\svgscale\undefined%
  \makeatother%
  \begin{picture}(1,0.3681547)%
    \put(0,0){\includegraphics[width=\unitlength,page=1]{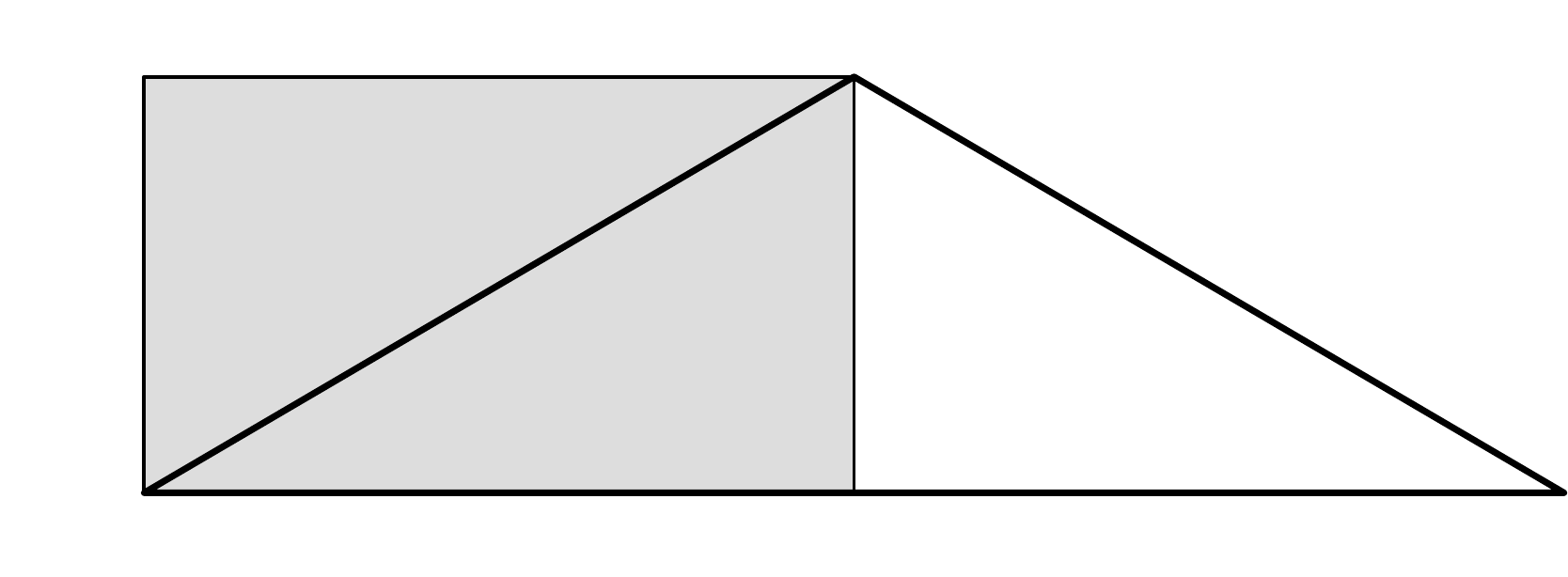}}%
    \put(0.26269082,0.00673395){\color[rgb]{0,0,0}\makebox(0,0)[lb]{\smash{$D$}}}%
    \put(0.26423694,0.33334396){\color[rgb]{0,0,0}\makebox(0,0)[lb]{\smash{$D$}}}%
    \put(0.5562299,0.1771392){\color[rgb]{0,0,0}\makebox(0,0)[lb]{\smash{$N$}}}%
    \put(-0.00380323,0.1771392){\color[rgb]{0,0,0}\makebox(0,0)[lb]{\smash{$N$}}}%
    \put(0.21693386,0.24102927){\color[rgb]{0,0,0}\makebox(0,0)[lb]{\smash{}}}%
    \put(0,0){\includegraphics[width=\unitlength,page=2]{disegno2.pdf}}%
  \end{picture}%
\endgroup%

\caption{In the Figure the letter $D$ and $N$ corresponds to Dirichlet boundary conditions and Neumann boundary conditions respectively. On such a rectangle a test function is provided by reflecting (anti-symmetrically) the eigenfunction of the triangle, from the light grey part to the dark grey part.}\label{fig2}
\end{figure}
For $\dfrac{\pi}{3} \le \alpha \le \pi$ we find that
\begin{equation}
\label{L1}
\lambda(\triangle_{\beta}) \ge \min \left\{\pi^2\left(1 + \frac{1}{d^2}\right), \frac{4 \pi^2}{d^2}\right\} = \pi^2\left(1 + \frac{1}{d^2}\right).
\end{equation}
Combining \eqref{T2} and \eqref{L1} we obtain
\begin{equation*}
\dfrac{T(\triangle_{\beta})\lambda(\triangle_{\beta}))}{|\triangle_{\beta}|} \ge \frac{\pi^2}{24},\, 0<\beta\le \dfrac{\pi}{3}.
\end{equation*}

Next we consider the case $0<\alpha \le\dfrac{\pi}{3}$ or $\pi/3\le \beta<\pi/2$. We have
\begin{equation}\label{T12}
|\triangle_{\beta}|=1/\tan(\alpha/2).
\end{equation}
 Let
\begin{equation*}
S(\rho,\alpha)=\{(r,\phi):0<r<\rho,\, -\alpha/2<\phi<\alpha/2\}
\end{equation*}
be the circular sector with radius $\rho$ and opening angle $\alpha$. Siudeja's Theorem 1.3 in \cite{BS} asserts that for $0<\beta\le \pi/3$, $\lambda(\triangle_{\pi/2-\alpha/2})\ge \lambda(S(\rho,\alpha))$,
where $d$ is such that $|\triangle_{\beta}|=|S(\rho,\alpha)|$. It follows that
\begin{equation}\label{T10}
\rho^2=2/(\alpha\tan(\alpha/2)).
\end{equation}
 Hence
\begin{equation*}
\lambda(\triangle_{\beta}) \geq 2^{-1}\alpha\tan(\alpha/2)j_{\pi/\alpha}^2.
\end{equation*}
where we have used that the first Dirichlet eigenvalue of a circular sector of opening angle $\beta$ and radius $\rho$ equals $ j^{2}_{\pi/\beta}\rho^{-2}$. See \cite{PSZ}.
Moreover by (1.2) and (4.3) for $k=1$ and $\nu=\pi/\alpha$ in \cite{QW}  we have
\begin{equation*}
j_{\pi/\alpha}^2 > \left(\dfrac{\pi}{\alpha}- \dfrac{a_1}{2^{1/3}} \left(\dfrac{\pi}{\alpha}\right)^{1/3}\right)^2, \quad -a_1 \ge \left(\dfrac{9 \pi}{8}\right)^{2/3},
\end{equation*}
where $a_1$ is the first negative zero of the Airy function. It follows that
\begin{equation}\label{T6}
j_{\pi/\alpha}^2\ge \frac{\pi^2}{\alpha^2}\bigg(1+C\bigg(\frac{\alpha}{\pi}\bigg)^{2/3}\bigg)^2\ge \frac{\pi^2}{\alpha^2}\big(1+C_1\alpha^{2/3}\big),
\end{equation}
where
\begin{equation}\label{T7}
C=(9\pi/8)^{2/3}2^{-1/3},\, C_1=(9/4)^{2/3}.
\end{equation}

The torsion function for $S(\rho,\alpha),\, \alpha<\pi/2,$ is given by (p.279 in \cite{TG}),
\begin{align*}
v_{S(\rho,\alpha)}(r,\phi)=&\frac{r^2}{4}\bigg(\frac{\cos(2\phi)}{\cos \alpha}-1\bigg)\nonumber \\ &+\frac{4\rho^2\alpha^2}{\pi^3}\sum_{n=1,3,5,...}(-1)^{(n+1)/2}\bigg(\frac{r}{\rho}\bigg)^{n\pi/\alpha}\cos\bigg(\frac{n\pi\phi}{\alpha}\bigg)
n^{-1}\bigg(n+\frac{2\alpha}{\pi}\bigg)^{-1}\bigg(n-\frac{2\alpha}{\pi}\bigg)^{-1}.
\end{align*}
By monotonicity of the torsion we obtain
\begin{align}\label{torsion}
T(\triangle_{\beta})& \geq T(S(\rho,\alpha))\nonumber \\ &=\int_{(0,d)}r\, dr\int_{(-\alpha/2,\alpha/2)}\, d\phi v_{S(d,\alpha)}(r,\phi)\nonumber \\ &
= \dfrac{d^4}{16} \bigg(\tan \alpha - \alpha - \dfrac{128 \alpha^4}{\pi^5} \sum\limits_{n=1,3,..} n^{-2}\bigg(n+\dfrac{2\alpha}{\pi}\bigg)^{-2}\bigg(n-\dfrac{2\alpha}{\pi}\bigg)^{-1}\bigg),
\end{align}
We have that for $0<\alpha\le \pi/3$, $(n+2\alpha/\pi)^{2}(n-2\alpha/\pi)\ge \frac{25}{27}n^3,\, n\in \N.$
This gives that
\begin{align}\label{T4}
T(\triangle_{\beta})&\ge \dfrac{\rho^4}{16}\bigg(\tan \alpha - \alpha - \dfrac{2^23^331\zeta(5)\alpha^4 d^4}{25\pi^5}\bigg)\nonumber \\ &
\ge \dfrac{\alpha^3\rho^4}{48}\big(1-C_2\alpha\big),
\end{align}
where
\begin{equation*}
C_2=\dfrac{2^23^431\zeta(5)}{5^2\pi^5}.
\end{equation*}
By \eqref{T6}, \eqref{torsion}, \eqref{T4}, and \eqref{T12} we obtain
\begin{equation}\label{finalestimate}
\dfrac{T(\triangle_{\beta})\lambda(\triangle_{\beta})}{|\triangle_{\beta}|}\ge \frac{\pi^2}{24}\big(1-C_2\alpha\big)\big(1+C_1\alpha^{2/3}\big).
\end{equation}
The right-hand side of \eqref{finalestimate} is greater or equal than $\frac{\pi^2}{24}$ for
\begin{equation}\label{T11}
C_1\ge C_1C_2\alpha+C_2\alpha^{1/3}.
\end{equation}
Inequality \eqref{T11} holds for all $\alpha\le 33/100.$

\subsection{Computer validation for the case $\beta\in (\pi/3,\beta_0)$ via interval arithmetic.}

We consider a triangle $\triangle^{\alpha}$ of height $1$ and opening angle $\alpha$, where $\alpha=\pi-2\beta$. Let
$$F(\alpha)=\frac{24}{\pi^2}\frac{\lambda(\triangle^\alpha)T(\triangle^\alpha)}{|\triangle^\alpha|}.$$
We wish to show that $F(\alpha)> { 1.01}$ in the range $0.33\le \alpha \le \pi/3$.

We present here a computer assisted proof of the result using Interval Arithmetic.

We once more use Siudeja's lower bound, comparing with the sector having the same opening angle and the same area (Theorem 1.3 of \cite{BS}), and get

\begin{equation}\label{La}
\lambda(\triangle^\alpha) \ge\tilde\lambda(\alpha)=\cos^2\left(\frac{\alpha}{2}\right)\left(\frac{\alpha}{\sin\alpha}\right)
\left(\frac{\pi}{\alpha}+C\left(\frac{\pi}{\alpha}\right)^\frac{1}{3}\right)^2,
\end{equation}
where $C$ is given by \eqref{T7}.

The area is given by
\begin{equation}\label{Ar}
|\triangle^\alpha|=\tan\left(\frac{\alpha}{2}\right),
\end{equation}

The monotonicity of $T$ with respect to inclusion allows us to estimate from below using the torsion of a tangent sector with same opening angle $\alpha$. We use \eqref{torsion} and find that
$$
T(\triangle^\alpha)\ge \frac{1}{16}\left(\tan\alpha-\alpha\right)-\frac{8}{\pi^5}\alpha^4\sum_{n=1,3,5,\dots} n^{-2}\left(n+\frac{2\alpha}{\pi}\right)^{-2}\left(n-\frac{2\alpha}{\pi}\right)^{-1}.
$$
In order to perform a numerical evaluation we truncate the series in the following way
$$
\begin{array}{l}
\displaystyle\sum_{n=1,3,5,\dots} n^{-2}\left(n+\frac{2\alpha}{\pi}\right)^{-2}\left(n-\frac{2\alpha}{\pi}\right)^{-1}=\sum_{n=0}^{\infty} (2n+1)^{-2}\left(2n+1+\frac{2\alpha}{\pi}\right)^{-2}\left(2n+1-\frac{2\alpha}{\pi}\right)^{-1}\\\\
\displaystyle\le\sum_{n=0}^{N} (2n+1)^{-2}\left(2n+1+\frac{2\alpha}{\pi}\right)^{-2}\left(2n+1-\frac{2\alpha}{\pi}\right)^{-1}+\frac{1}{2^5}\sum_{n=N+1}^{\infty}\frac{1}{n^5}\\\\
\displaystyle\le\sum_{n=0}^{N} (2n+1)^{-2}\left(2n+1+\frac{2\alpha}{\pi}\right)^{-2}\left(2n+1-\frac{2\alpha}{\pi}\right)^{-1}+\frac{1}{2^7N^4}.
\end{array}$$
It follows that
\begin{align}
T(\triangle^\alpha)&\ge \tilde{T}(\triangle^\alpha)\nonumber \\ &= \frac{1}{16}\left(\tan\alpha-\alpha\right)-\frac{8}{\pi^5}\alpha^4\left(
\sum_{n=0}^{10} (2n+1)^{-2}\left(2n+1+\frac{2\alpha}{\pi}\right)^{-2}\left(2n+1-\frac{2\alpha}{\pi}\right)^{-1}+\frac{1}{2^7\cdot 10^4}
\right).
\end{align}

Therefore
$$F(\alpha)\ge G(\alpha)=\frac{24}{\pi^2}\frac{\tilde{\lambda}(\triangle^\alpha)\tilde{T}(\triangle^\alpha)}{|\triangle^\alpha|}$$

At this point we can prove that $G(\alpha)>1.01$ for all values $0.33\le\alpha\le \pi/3$ by using Interval Arithmetic. There are many softwares and libraries which can be employed for this purpose. We selected \emph{Octave}\footnote{John W. Eaton, David Bateman, Søren Hauberg, Rik Wehbring (2018).
GNU Octave version 4.4.1 manual: A high-level interactive language for
numerical computations.
URL https://www.gnu.org/software/octave/doc/v4.4.1/} (A free software that runs on GNU/Linux, macOS, BSD, and Windows) which provides a specific package called \emph{Interval}.\footnote{Oliver Heimlich, GNU Octave Interval Package, https://octave.sourceforge.io/interval/, version 3.2.0, 2018-07-01. The interval package is a collection of functions for interval arithmetic.  It is developed at Octave Forge, a sibling of the GNU Octave project.}

We covered the interval $\left[\frac{33}{100},\frac{\pi}{3}\right]$ by a collection of 1001 intervals $I_n$ with $n=0,\dots, 1000$, so that
$$I_n=\left[\frac{33}{100} + \frac{(n-1)}{1000}\left(\frac{\pi}{3}-\frac{33}{100}\right),\frac{33}{100} + \frac{(n+1)}{1000}\left(\frac{\pi}{3}-\frac{33}{100}\right)\right].$$
We observe that the intersection of consecutive intervals is intentionally non empty.
Using the \emph{Interval} package,  we designed a code that for $n$ going from $0$ to $10^3$ 
provides upper and lower bounds for $G(I_n)$ in terms of floating point numbers. This is performed in an automated way by standard and reliable algorithms. We established that the inequality
$F(\alpha)>1.01$ holds true on the whole interval $\left[\frac{33}{100},\frac{\pi}{3}\right]$ by verifying it on $I_n$ for all $n=0,\dots,10^3$. 
\hspace*{\fill }$\square $
\vspace{5mm}

For completeness we include the code below.

 \lstset{language=Octave,backgroundcolor=\color{mygray},frame=single,numbers=left,breaklines=true,keywordstyle=\bfseries\color{red}, identifierstyle=\color{blue},stringstyle=\color{violet}}
 \begin{lstlisting}
pkg load interval # load the package Interval
output_precision (6) # number of digits displayed
C=(9*pi/8)^(2./3)*2^(-1./3);

function K=G(x) # this is the definition of the function G(\alpha)
  Sum=0;
  for n = 0:10
    Sum = Sum + (2*n+1)^(-2)*(2*n+1+2*x/pi)^(-2)*(2*n+1-2*x/pi)^(-1);
  endfor
  K=(24./pi^2)*(1./16*(tan(x)-x)-8*x^4/pi^5*(Sum+1./(2.^7*10.^4)))*((cos(x/2))^2*(x/sin(x))*(pi/x+((9*pi/8)^(2./3)*2^(-1./3))*(pi/x)^(1/3))^2)/tan(x/2);
endfunction
control="OK"; # the variable control is set to "OK"
N=1000; # Number of intervals
Delta=(pi/3-0.33)/N; # 2Delta is the size of each interval
for n = 0:N
  n #print the value of n
  a=0.33+n*Delta;
  I=midrad(a,Delta) # I = interval with center in a and radius Delta
  J=G(I) # J is an interval which includes the image of I
  if(J>1.01) # check that (min J) > 1.01
  "so far inequality holds" # tell that everything is working fine
  elseif
  control="failure" # the variable control is set to "failure"
  break # in case of failure the cycle breaks
  endif
endfor
\end{lstlisting}

{\bf Acknowledgments.} The authors wish to thank Prof. Gerardo Toraldo for helpful discussions on Interval Arithmetic. The authors acknowledge support by the London Mathematical Society, Scheme 4 Grant 41636. MvdB was also supported by The Leverhulme Trust through Emeritus Fellowship EM-2018-011-9, and by an INDAM-GNAMPA Grant for visiting professors.

\end{document}